\documentclass[a4paper,leqno]{amsart}

\usepackage{latexsym}
\usepackage[english]{babel}
\usepackage{fancyhdr}
\usepackage[mathscr]{eucal}
\usepackage{amsmath}
\usepackage{mathrsfs}
\usepackage{amsfonts}
\usepackage{amssymb}
\usepackage{amscd}
\usepackage{bbm}
\usepackage{graphicx}
\usepackage{subcaption}
\usepackage{graphics}
\usepackage{latexsym}
\usepackage{color}

\usepackage{pifont}
\usepackage{booktabs} 

\newcommand{\ud}{\mathrm{d}}

\newcommand{\cH}{\mathcal{H}}
\newcommand{\ran}{\mathrm{ran}}

\newcommand{\N}{\mathbb N}

\theoremstyle{plain}
\newtheorem{theorem}{Theorem}[section]

\newtheorem{corollary}[theorem]{Corollary}
\newtheorem{proposition}[theorem]{Proposition}

\theoremstyle{definition}
\newtheorem{definition}[theorem]{Definition}

\newtheorem{example}[theorem]{Example}

\numberwithin{equation}{section}



\title[On Krylov solvability of Friedrichs systems]{Open problems and perspectives on solving Friedrichs systems by Krylov approximation}
\author[N.~A.~Caruso]{No\`e Angelo Caruso}
\address[N.~A.~Caruso]{Dipartimento di Scienze Umane, Universit\`{a} degli Studi ``Link Campus University'' \\ Via del Casale di San Pio V, 44 \\ 00165 Rome (Italy)}
\email{n.caruso@unilink.it} 

\author[A.~Michelangeli]{Alessandro Michelangeli} 
\address[A.~Michelangeli]{Department of Mathematics and Natural Sciences, Prince Mohammad Bin Fahd University \\ Al Khobar 31952 (Saudi Arabia) \\
\and
Hausdorff Center for Mathematics, University of Bonn \\ Endenicher Allee 60 \\ D-53115 Bonn (Germany) \\ 
\and 
TQT Trieste Institute for Theoretical Quantum Technologies, Trieste (Italy)}
\email{amichelangeli@pmu.edu.sa}

\begin{document}

\begin{abstract}
 We set up, at the abstract Hilbert space setting, the general question on when an inverse linear problem induced by an operator of Friedrichs type admits solutions belonging to (the closure of) the Krylov subspace associated to such operator. Such Krylov solvability of abstract Friedrichs systems allows to predict when, for concrete differential inverse problems, truncation algorithms can or cannot reproduce the exact solutions in terms of approximants from the Krylov subspace.
 \end{abstract}



\maketitle


\section{Introduction}\label{sec:intro-main}

At the intersection of the two general fields of numerical solutions to partial differential equations, and of iterative approximations of inverse linear problems, an interesting class of questions has been recently emerging, with particularly rich developments and impact both in an abstract setting and in concrete applications. In this note we would like to briefly review this subject and, in a sense, to advocate for it.

In abstract terms, an \emph{inverse linear problem} in some linear space $\cH$ consists of the following: an element $g\in\cH$ (the `datum') and a linear map $A$ on $\cH$ are given, and one is to find, if any, those element(s) $f\in\cH$ such that $Af=g$.

This clearly encompasses a huge variety of concrete classes of inverse problems, including those where $A$ is a differential operator. On the one hand, it is customary to study each such class developing ad hoc methods and tools. This falls within the realm of numerical analysis, and includes the issue of reproducing the `exact' solution(s) $f$ to $Af=g$ as a limit of approximants $f_n\to f$ obtained by solving (easier) truncated problems $A_nf_n=g_n$, typically in an iterative way.

Another perspective is to disregard the specific features and setting of inverse linear problems and to investigate them at a sufficiently high level of abstract generality, namely to study the solvability of $Af=g$ and the approximability $f_n\to f$ depending on: the finite or infinite dimensionality of the vector space $\cH$, the additional analytic or topological structure of $\cH$ (e.g., a Banach space, a Hilbert space), the boundedness or unboundedness of the operator $A$ on $\cH$ and further features of $A$ such as its compactness or normality or self-adjointness and the like, and the class of data $g$ with respect to the given operator $A$.

Obviously, there is a natural interplay between the two perspectives above. For instance, the actual convergence $f_n\to f$ from solutions to suitably truncated inverse problems obtained in a concrete setting, as well as the convergence rate and the quantitative control (say, the norm) of the displacement $f-f_n$, are often established by means of abstract functional-analytic and operator theoretic methods. We mention, in this respect, excellent classical references such as \cite{kranoselski-1972-approxsoll,Nevanlinna-1993-converg-iterat-book,Vainikko-1993,Hanke-ConjGrad-1995,Engl-Hanke-Neubauer-1996,Chen3-1997,Hansen-Illposed-1998,Liesen-Strakos-2003,Saad-2003_IterativeMethods,Ern-Guermond_book_FiniteElements,Atkinson-Han-TheoNumAnal2009,Sullivan-UncertQuant,Quarteroni-book_NumModelsDiffProb}, and also, to quickly grasp the flavour of the abstract point of view, \cite[Sections 1.1 and 1.2, and Appendix A]{CM-KrylovMonograph-2021}.

Our interest here concerns the abstract analysis of a class of approximation schemes for inverse linear problems known as \emph{Krylov subspace methods}. Their efficiency and popularity places them among the `Top 10 Algorithms' of the 20th century \cite{Dongarra-Sullivan-Best10-2000,Cipra-SIAM-News}.

In a nutshell (Section \ref{sec:crashcourseKrylov} provides a concise and somewhat more detailed `crash course' on the subject), one assumes that $\cH$ is a Hilbert space (clearly, the richest scenario of problems and results occurs when $\dim\cH=\infty$), and given the operator $A$ acting on $\cH$ and the vector $g\in\mathrm{ran}A$ one suitably truncates the inverse linear problem $Af=g$ to the finite-dimensional subspaces $\mathrm{span}\{g,Ag,A^2g,\dots,A^{n-1}g\}$ for larger and larger integers $n$. If $f_n$ is a solution of the $n$-truncated problem, then under favourable circumstances a solution $f$ to $Af=g$ is approximated in norm by the $f_n$'s, i.e., $\|f-f_n\|_{\cH}\to 0$ as $n\to\infty$, which means, in particular, that $f$ belongs to the closure of the subspace generated by $g,Ag,A^2g,\dots$

Thus, Krylov subspace methods are (typically iterative) algorithms where a solution $f$ to an inverse linear problem $Af=g$ is obtained by closer and closer approximants that are expressed by convenient finite linear combinations of the vectors $g,Ag,A^2g,\dots$. The subspace
\begin{equation}\label{eq:Krylovsubspace}
 \mathcal{K}(A,g)\,:=\,\mathrm{span}\big\{ g, Ag, A^2g, \dots \big\}\,\subset\,\cH
\end{equation}
is called in this context the \emph{Krylov subspace} associated with $A$ and $g$ (in other contexts it is also called the \emph{cyclic subspace} for $A$ and $g$). When the inverse problem $Af=g$ admits a solution $f\in\overline{\mathcal{K}(A,g)}$, the problem is said to be \emph{Krylov-solvable}.

In a sense, for inverse problems that are Krylov-solvable, the advantage is that the trial space for the approximants is spanned by vectors (the $A^n g$'s) that are explicit and typically easy to calculate in the applications. This makes Krylov subspace methods particularly versatile.

The study of the Krylov solvability of inverse linear problems at an abstract Hilbert space setting, and of the asymptotic convergence of the Krylov approximants, is a fascinating subject that only very recently has been set up in a thorough and systematic way. We refer on this to the monograph \cite{CM-KrylovMonograph-2021} and to the `crash course' of Section \ref{sec:crashcourseKrylov}.

As a matter of fact, the theoretical apparatus developed so far in the subject of Krylov solvability is sufficiently articulated to make it natural to specialise it for distinguished classes of inverse problems. The one we consider here is the type of inverse problems customarily referred to as \emph{Friedrichs systems}.

Friedrichs systems are a wide class of differential linear inverse problems essentially characterised by the fact that the sum of the differential operator $A$ and its formal adjoint is bounded and coercive. They have recently attained particular relevance in numerical analysis \cite{HMSW,MJensen,Ern-Guermond_book_FiniteElements,EGsemel} as a convenient unified framework for numerical solutions to partial differential equations of different types, including elliptic, parabolic, and hyperbolic.

Starting from concrete, differential formulations of Friedrichs-type inverse problems, it is possible to re-formulate them in an abstract Hilbert space setting, so as to exploit general operator-theoretic methods and then apply the results to concrete versions \cite{Ern-Guermond_book_FiniteElements,EGC,ABcpde,ABCE,Antonic-Erceg-Mich-2017,EM-FriedrichsDelta2017,Erceg-Soni-2022}. Section \ref{sec:overviewFriedrichs} revisits the basic aspects of this construction.

This naturally leads to the notion of a \emph{joint pair of abstract Friedrichs operators} in a Hilbert space $\cH$, namely a pair $(A_0,\widetilde{A}_0)$ of two densely defined and closed linear operators $A_0:\mathcal{W}_0\to\cH$ and $\widetilde{A}_0:\mathcal{W}_0\to\cH$ such that
\begin{equation}\label{eq:defFriedrichs}
\begin{split}
 & A_0\,\subset\,(\widetilde{A}_0)^*\,=:A_1\,,\qquad \widetilde{A}_0\,\subset\,(A_0)^*\,=:\widetilde{A}_1\,, \\
 & \textrm{$A_1+\widetilde{A}_1$ is a bounded self-adjoint operator in $\cH$} \\
 & \textrm{with strictly positive bottom}\,.
\end{split}
\end{equation}

At this abstract level, $A_0$ and $A_1$ model what in a typical concrete example of inverse problem of Friedrichs type are, respectively, the minimal and the maximal realisation of a differential operator on a bounded domain $\Omega\subset\mathbb{R}^n$. Any closed operator $A$ with $A_0\subset A\subset A_1$ models instead a specific realisation of the differential operator, identified by specific boundary conditions at the boundary $\partial\Omega$.

At the overlap of the two seemingly distant subjects of abstract Krylov solvability and abstract Friedrichs systems we finally make the following observation: \emph{investigating in the abstract Hilbert space setting the Krylov solvability of an inverse problem $Af=g$ governed by an operator $A$ of Friedrichs type would yield extremely valuable a priori information on the possibility of approximating numerically the solution to a concrete differential inverse problem within the vast Friedrichs class}.

This note is conceived to serve both as a self-contained introduction to the set of problems we outlined above, including the current state of the art and the recent achievements, and as a survey of relevant questions on top of the research agenda in this subject.

As mentioned already, Sections \ref{sec:crashcourseKrylov} and \ref{sec:overviewFriedrichs} provide background materials. In Section \ref{sec:nextresultsandproblems} we outline a few of the most relevant and natural questions that are worth addressing.

\section[Short `crash course' on Krylov solvability]{A very short `crash course' on Krylov solvability in abstract Hilbert space}\label{sec:crashcourseKrylov}

For the concise review materials of this Section we refer to the corpus of recent works \cite{CMN-truncation-2018,CMN-2018_Krylov-solvability-bdd,CM-Nemi-unbdd-2019,CM-2019_ubddKrylov,CM-KrylovPerturbations-2020,CM-KrylovMonograph-2021,C-KrylovNormal-2022}.

\subsection{Krylov approximation methods as truncation algorithms}\label{subsec:Krylov-truncation}~

We begin by presenting Krylov subspace methods as a family of approximation methods based on suitable compressions, or `truncations', of the operator $A$ and the vector $g$ of the inverse linear problem. Certainly one may think of other ways to represent various Krylov based methods, such as the more well-known iterative or variational formulations, yet these viewpoints have the drawback that they somewhat mask the relevance of the Krylov subspace and the concept of `Krylov solvability'.

Let us recall from the Introduction that the inverse linear problem, in the complex Hilbert space $\cH$, 
\begin{equation}\label{eq:invlin}
Af = g
\end{equation}
in the unknown $f\in\cH$ and for given densely defined (and possibly closed) linear operator $A$, and $g \in \ran\,A\subset \cH$ -- which ensures the existence of solution(s) -- is said to be \emph{Krylov-solvable} if there is a solution $f$ belonging to the closure in $\cH$ of the Krylov subspace
\begin{equation}\label{eq:Krylovsubspace2}
 \mathcal{K}(A,g)\,:=\,\mathrm{span}\big\{ g, Ag, A^2g, \dots \big\}\,\subset\,\cH
\end{equation}
associated to $A$ and $g$. This expresses the notion of \emph{Krylov solvability} of \eqref{eq:invlin} and $f$ in this case is referred to as a \emph{Krylov solution}.


Krylov methods work by constructing a sequence of approximants $(f_n)_{n \in \N}$ of $f$ that are solutions to a \emph{truncated} problem obtained by `projecting' the original \eqref{eq:invlin} onto an $n$-dimensional subspace built with a finite amount of vectors of the type $A^k g$.

In practice one considers the subspace
\begin{equation}
 \mathcal{K}_n(A,g)\,:=\,\mathrm{span}\big\{ g, Ag, A^2g, \dots , A^{n-1}g \big\}\,,\qquad \dim\mathcal{K}_n(A,g)\leqslant n\,,
\end{equation}
as well as another, possibly distinct, $n$-dimensional subspace $\mathcal{J}_n(A,g)$ consisting of finite linear combinations of an amount of the $A^kg$'s. Correspondingly, consider the orthogonal projections $P_n:\cH \to \cH$ and $Q_n:\cH \to \cH$, respectively, onto $\mathcal{K}_n(A,g)$ and $\mathcal{J}_n(A,g)$. It is customary to refer to $\mathcal{K}_n(A,g)=\mathrm{ran}\,P_n$ and $\mathcal{J}_n(A,g)=\mathrm{ran}\,Q_n$ as, respectively, the \emph{solution subspace} and the \emph{trial subspace}.

The $n$-truncated inverse linear problem is then formulated as 
\begin{equation}\label{eq:invlin_truncated}
 \begin{split}
  A_nf_n\,&=\,g_n\,, \\
  \textrm{where} & \; A_n\,:=\,Q_n A P_n\,,\quad g_n \,:=\, Q_n g\,,
 \end{split}
\end{equation}
 and the $n\times n$ matrix $A_n$ is the \emph{compression} of $A$ to the considered subspaces. It need \emph{not} happen that $f_n = P_n f$, where $f$ is a solution to the `full' problem \eqref{eq:invlin}. The finite-dimensional inverse linear problems \eqref{eq:invlin_truncated} form the basis of many different approximation schemes, in particular the famous Petrov-Galerkin methods.

 As we increase the dimension of the subspaces, i.e., as $n$ becomes large, ideally an approximation scheme produces solutions to the truncated problem $(f_n)_{n \in \N}$ that tend to an exact solution of the original inverse linear problem \eqref{eq:invlin}. The convergence $f_n\to f$ in the norm of $\cH$ is tantamount of the Krylov solvability of \eqref{eq:invlin} defined above.

 We remark that the choice of $Q_n$ (and hence of $\mathcal{J}_n(A,g)$) may vary depending on the precise method used. For example, the celebrated GMRES method uses the choice of $Q_n$ as the orthogonal projection operator onto the subspace $A\mathcal{K}_n(A,g)$. On the other hand, for the class of problems for which $A$ is positive self-adjoint, the famous conjugate gradient (CG) method uses the choice $Q_n = P_n$.

 Beside the above-mentioned recent discussion \cite{CMN-truncation-2018}, an excellent overview and general results associated with the formulation and solution to \eqref{eq:invlin_truncated} may be found in several excellent articles and monographs such as \cite{Chen3-1997,kranoselski-1972-approxsoll,Nevanlinna-1993-converg-iterat-book,Vainikko-1993,vainikko1981,vainikko1974}.

\subsection{A glance at Krylov solvability for bounded inverse problems}\label{subsec:Krylov_solvability}~

 Interesting structural properties emerge in the discussion of Krylov solvability of inverse problems \eqref{eq:invlin} with $A \in \mathcal{B}(\cH)$ (for a more systematic discussion, see \cite{CMN-2018_Krylov-solvability-bdd}).
 
 Let us start with the notion of \emph{Krylov reducibility}. The following inclusions are easily seen:
 \begin{equation}
  \begin{split}
   A\,\overline{\mathcal{K}(A,g)}\,& \subset\, \overline{\mathcal{K}(A,g)}\,, \\
   A^*\mathcal{K}(A,g)^\perp \,& \subset\, \mathcal{K}(A,g)^\perp\,,
  \end{split}
 \end{equation}
  which express the invariance of $\overline{\mathcal{K}(A,g)}$, respectively, of $\mathcal{K}(A,g)^\perp$, under the action of $A$, respectively, of $A^*$. One then says that $A$ is \emph{$\mathcal{K}(A,g)$-reduced} if in addition $\mathcal{K}(A,g)^\perp$ is invariant also under the action of $A$, i.e.,
  \begin{equation}
   A\,\mathcal{K}(A,g)^\perp \,\subset\, \mathcal{K}(A,g)^\perp\,.
 \end{equation}
  Here \emph{Krylov reducibility} indicates that $A$ is reduced (in the standard operator sense, e.g., \cite[Section 1.4]{schmu_unbdd_sa}) into two blocks with respect to the Hilbert space direct sum decomposition
  \begin{equation}
   \cH \,=\,\overline{\mathcal{K}(A,g)}\,\oplus\,\mathcal{K}(A,g)^\perp\,.
  \end{equation}

  An evident sufficient condition for Krylov reducibility is self-adjointness: a bounded self-adjoint operator $A$ is \emph{always} $\mathcal{K}(A,g)$-reduced for any $g \in \cH$. More generally, bounded normal operators $A$ are $\mathcal{K}(A,g)$-reduced if and only if $A^*g\in\overline{\mathcal{K}(A,g)}$ \cite[Proposition 2.1]{CMN-2018_Krylov-solvability-bdd}. Thus, for unitary operators Krylov reducibility is equivalent to Krylov solvability \cite[Remark 3.6]{CMN-2018_Krylov-solvability-bdd}.
 
  It turns out that the Krylov reducibility of $A$ is intimately connected to the Krylov solvability of the associated  inverse linear problem \eqref{eq:invlin}.

\begin{proposition}\cite[Proposition~3.9]{CMN-2018_Krylov-solvability-bdd}\label{prop:Kreduced_Ksolvable}
Let $A \in \mathcal{B}(\cH)$ and let $g \in \ran A$. If A is $\mathcal{K}(A, g)$-reduced, then there exists a Krylov solution to the problem $A f = g$. For example, if $f \in \cH$ satisfies $Af = g$ and $P_\mathcal{K}$ is the orthogonal projection onto
$\overline{\mathcal{K}(A, g)}$, then $f^\circ := P_\mathcal{K} f$ is a Krylov solution
\end{proposition}

The immediate corollary of the above proposition is:

\begin{corollary}\label{cor:bddself-adjoint_Ksolvable}
Let $A \in \mathcal{B}(\cH)$ be self-adjoint and let $g \in \ran A$. Then there exists a solution to the inverse linear problem $A f = g$ in the space $\overline{\mathcal{K}(A,g)}$.
\end{corollary}

Therefore, the bounded-self adjoint inverse problems are \emph{always} Krylov-solvable.

 As discovered in \cite{CMN-2018_Krylov-solvability-bdd}, the following subspace is inherently linked to the Krylov solvability of an inverse linear problem.

\begin{definition}\label{def:Kintersection}
The \emph{Krylov intersection} associated to given $A \in \mathcal{B}(\cH)$ and $g \in \cH$ is the subspace
\begin{equation}\label{eq:def_Kint}
\mathcal{I}(A,g) \,:=\, \overline{\mathcal{K}(A,g)} \cap A (\mathcal{K}(A,g)^\perp)\,.
\end{equation}
\end{definition}

 The possible triviality of the Krylov intersection is another structural property relevant for the Krylov solvability.

\begin{proposition} \cite[Proposition~3.4]{CMN-2018_Krylov-solvability-bdd},  \cite[Proposition~2.4]{CM-2019_ubddKrylov}\label{prop:trivialKint_Ksolv}
Let $A \in \mathcal{B}(\cH)$ and let $g \in \ran A$. If $\mathcal{I}(A,g) = \{0\}$, then there exists a Krylov solution $f \in \overline{\mathcal{K}(A,g)}$ to the inverse linear problem $A f = g$. If, additionally, $A^{-1} \in \mathcal{B}(\cH)$, then $\mathcal{I}(A,g) = \{0\}$ if and only if the solution $f$ to $Af = g$ is in $\overline{\mathcal{K}(A,g)}$.
\end{proposition}

 Observe that, obviously, the $\mathcal{K}(A,g)$-reducibility of an operator $A$ is just a special case of the triviality of the Krylov intersection $\mathcal{I}(A,g)$, which indicates that the latter subspace is more informative concerning the possibility that an inverse problem be Krylov-solvable.

 We refer to \cite{CMN-2018_Krylov-solvability-bdd} for an expanded discussion of classes of inverse linear problems that are Krylov-reduced, or have trivial Krylov intersection, and in general are Krylov-solvable, including when the Krylov solution is unique. To that, we add the following interesting recent result: 
 
 \begin{proposition} \cite[Propositions~3.4 and 3.9]{C-KrylovNormal-2022}\label{prop:Acompactnormal_Ksolvable}
Let $ A \in \mathcal{B}(\cH)$ be a compact normal operator and let $g \in \cH$. Then $A$ is $\mathcal{K}(A, g)$-reduced. If, additionally, $g \in \ran A$, then $A f = g$ has a unique Krylov solution that is also a minimal norm solution to the inverse linear problem.
\end{proposition}

 \subsection{Krylov solvability and unbounded self-adjoint operators}\label{subsec:Krylovsolvability_unbddselfadj}~

 As a matter of fact, the notion of Krylov reducibility and Krylov intersection, in connection to the Krylov solvability of an inverse linear problem, have a counterpart in the case of \emph{unbounded} $A$, with analogous results to those outlined above: this has been discussed in \cite{CM-2019_ubddKrylov}.

 The unavoidable domain issues that arise in the unbounded case require suitable adaptations. For one thing, the very notion \eqref{eq:Krylovsubspace2} of the Krylov subspace $\mathcal{K}(A,g)$ only makes sense if $g$ belongs to the domain of \emph{every} power of $A$, an inescapable minimal condition when $A$ is unbounded. When this happens, the generalised notion of Krylov intersection reads
 \begin{equation}\label{eq:def_KInt_unbdd}
\mathcal{I}(A,g) := \overline{\mathcal{K}(A,g)} \cap A( \mathcal{K}(A,g)^\perp \cap \mathrm{dom}\, A )\,,
\end{equation}
 and a discussion of its properties and relevance for the Krylov solvability of $Af=g$ may be found in \cite[Sections 5 and 6]{CM-2019_ubddKrylov}.

 In fact, distinguished classes of data $g$ turn out to play a role in the analysis of the Krylov solvability of $Af=g$, and the following classical definition collects some of them.

 \begin{definition}\label{def:bdd_a_qa_class}
Let $A$ be a closed operator acting in a complex Hilbert space $\cH$ and let $g\in\cH$.
 \begin{enumerate}
  \item[(i)] $g$ belongs to the class $C^\infty(A)$ of \emph{smooth vectors} for $A$ if
  \[
   g\,\in\,\mathrm{dom}\,A^n\qquad \forall n\in\mathbb{N}\,.
  \]
  \item[(ii)] $g$ belongs to the class $\mathcal{D}^b(A)$ of \emph{bounded vectors} for $A$ if $g$ is smooth for $A$ and there is a constant $B_g>0$ such that
  \[
   \| A^n g\|_{\cH}\,\leqslant\, (B_g)^n\qquad \forall n\in\mathbb{N}\,.
  \]
  \item[(iii)] $g$ belongs to the class $\mathcal{D}^a(A)$ of \emph{analytic vectors} for $A$ if $g$ is smooth for $A$ and there is a constant $C_g>0$ such that
  \[
   \| A^n g\|_{\cH}\,\leqslant\, (C_g)^n\, n!\qquad \forall n\in\mathbb{N}\,.
  \]
  \item[(iv)] $g$ belongs to the class $\mathcal{D}^{qa}(A)$ of \emph{quasi-analytic vectors} for $A$ if $g$ is smooth for $A$ and 
  \[
   \sum_{n = 0}^\infty \Vert A^n g \Vert^{-\frac{1}{n}}_\cH \,=\, +\infty\,.
   \]
 \end{enumerate}
\end{definition}

 Clearly, 
 \begin{equation}
  \mathcal{D}^b(A) \,\subset\, \mathcal{D}^a(A) \,\subset\,  \mathcal{D}^{qa}(A) \,\subset\,  C^\infty(A) \,.
 \end{equation}

 The self-adjoint case exhibits in this respect distinguished properties. For one thing, as well known (see, e.g., \cite[Lemma 7.13]{schmu_unbdd_sa}), if $A$ is self-adjoint, then $\mathcal{D}^b(A)$ is dense in $\cH$, and so too are, consequently, all the other classes above. Moreover, a noticeable conclusion can be drawn concerning the Krylov solvability of possibly \emph{unbounded} self-adjoint inverse problems.

\begin{theorem} \cite[Theorems~4.1 and 4.2]{CM-2019_ubddKrylov}\label{th:SelfSkew-adjoint_Ksolvability}
Let $A$ be a self-adjoint or skew-adjoint operator on a Hilbert space $\cH$, and let $g \in \ran A \cap \mathcal{D}^{qa}(A)$. Then there exists a unique Krylov solution to the inverse linear problem $Af = g$. 
\end{theorem}

  It is worth commenting that Theorem \ref{th:SelfSkew-adjoint_Ksolvability} relies on an ingenious and somewhat laborious analysis carried out in the framework of the convergence of the conjugate gradient algorithm on an infinite-dimensional Hilbert space. This was initially developed in the \emph{bounded} case \cite{Nemirovskiy-Polyak-1985,Nemirovskiy-Polyak-1985-II} and recently re-designed in the \emph{unbounded} setting \cite[Theorem~2.4 and Corollary~2.5]{CM-Nemi-unbdd-2019}. The full proof obtained in \cite{CM-Nemi-unbdd-2019} in fact reveals that \emph{the sequence of approximants from the conjugate-gradients algorithm converges to the minimal norm solution of the inverse problem \eqref{eq:invlin}} when $A$ is positive self-adjoint and $g \in \mathcal{D}^{qa}(A) \cap \ran A$. In other words: the truncation scheme and truncated inverse problems \eqref{eq:invlin_truncated} defined by setting the trial space $Q_n = P_n$ for $P_n$ the orthogonal projection onto $\mathcal{K}_n(A,g)$, produce a well-defined sequence of approximants $(f_n)_{n \in \N}$ and $f_n \xrightarrow{n\to\infty} f^\circ$ in $\cH$, where $f^\circ$ is the minimal norm solution to \eqref{eq:invlin}.

\subsection{Krylov subspaces and solvability under perturbations}\label{subsec:KrylovPerturb}~

A further approach to the question of determining the Krylov solvability (or lack-of) of an inverse linear problem $Af=g$ is to analyse possibly `easier' \emph{perturbed} problems $\widetilde{A}\widetilde{f}=\widetilde{g}$, where in some suitable sense $\widetilde{A}\sim A$ and $\widetilde{g}\sim g$.

Even elementary toy problems show that the phenomenology may be pretty diverse.


\begin{example} \cite[Example~3.3]{CM-KrylovPerturbations-2020}\label{eg:Krylovperturb_lossgain}
Let $R : \ell^2(\mathbb{Z}) \to \ell^2(\mathbb{Z})$ be the right shift operator $e_n \mapsto e_{n + 1}$ (where $e_n$ is the $n$-th canonical vector of $\ell^2(\mathbb{Z})$). $R$ is unitary with $R^* = R^{-1} = L$, the left-shift operator $e_{n} \mapsto e_{n - 1}$. $R$ admits a dense set $\mathcal{C} \subset \ell^2(\mathbb{Z})$ of cyclic vectors (i.e., vectors $g \in \cH$ such that $\mathcal{K}(R,g)$ is dense in $\ell^2(\mathbb{Z})$) and a dense subset $\mathcal{N} \subset \ell^2(\mathbb{Z})$ of all finite linear combinations of the canonical vectors. Clearly all vectors in $\mathcal{N}$ are non-cyclic for the operator $R$, and moreover the solution $f$ to $R f = g$ for any $g \in \mathcal{N}$ does not belong to $\overline{\mathcal{K}(R,g)}$.
\begin{itemize}
	\item[(i)] (Loss of Krylov solvability.) For a datum $g \in \mathcal{N}$, the inverse problem $Rf = g$ admits a unique solution $f$ that is not in $\overline{\mathcal{K}(R,g)}$. Yet, by the density of $\mathcal{C}$ in $\ell^2(\mathbb{Z})$, there exists a sequence $(g_{m})_{m \in \N} \subset \mathcal{C}$ such that $\Vert g - g_m \Vert_{\ell^2} \xrightarrow{m \to \infty} 0$ and moreover for each $m \in \N$ the inverse linear problem $R f_m = g_m$ is Krylov-solvable with the unique solution $f_m = R^{-1}g_m$. Here we see that Krylov solvability is lost in the limit $m \to \infty$, although $\Vert f_m - f \Vert_{\ell^2} \xrightarrow{m \to \infty} 0$.
	\item[(ii)] (Gain of Krylov solvability.) For a datum $g \in \mathcal{C}$ the inverse linear problem $Rf = g$ is Krylov-solvable as $\overline{\mathcal{K}(R,g)} = \ell^2(\mathbb{Z})$. Yet, by density there exists a sequence $(g_m)_{m \in \N} \subset \mathcal{N}$ such that $\Vert g_m - g \Vert_{\ell^2} \xrightarrow{m \to \infty} 0$, and each perturbed problem $R f_m = g_m$ is not Krylov-solvable. Therefore, Krylov solvability is absent at the perturbed level but appears only in the limit, still with the approximation $\Vert f_m - f \Vert_{\ell^2} \xrightarrow{m \to \infty} 0$.
\end{itemize}
\end{example}

In a perturbative perspective, it is relevant to compare the possibly ``close'' Krylov subspaces $\mathcal{K}(A,g)$ and $\mathcal{K}(\widetilde{A},\widetilde{g})$. A promising tool is the \emph{weak gap metric} between subspaces of  \emph{separable} Hilbert space, inspired to the standard \emph{norm gap distance} (namely, the Hausdorff distance) between subspaces (see, e.g., \cite[Chapt.~4, \S 2]{Kato-perturbation}). Owing to the separability, and also reflexivity, of $\cH$, the Hilbert space weak topology restricted to the (norm-closed) unit ball $B_\cH$ of $\cH$ makes the ball a compact and complete metric space $(B_\cH,\rho_w)$, and the metric $\rho_w$ is induced by a norm, namely $\rho_w(x,y)=\|x-y\|_w$ for some explicitly constructed norm in $B_\cH$. This allows to define, for any pair of norm-closed linear subspaces $U$ and $V$ of $\cH$,
\begin{equation}\label{eq:defweakgap}
\left\{ \begin{aligned} 
d_w(B_U,B_V) & := \sup_{u \in B_U} \inf_{v \in B_V} \rho_w(u,v)\,, \\
\widehat{d}_w(U,V) & := \max\{ d_w(B_U,B_V), d_w(B_V,B_U) \} \,, \end{aligned} \right. 
\end{equation}
where $B_U = B_\cH \cap U$ and $B_V = B_\cH \cap V$. As a matter of fact \cite[Lemmas 6.1 and 6.2]{CM-KrylovPerturbations-2020}, $\widehat{d}_w$ is indeed a metric for the norm-closed subspaces of $\cH$, albeit in general it is \emph{not} complete.

In the context of Krylov subspaces, the distance $\widehat{d}_w$ is effective in a multiplicity of instances, for example in monitoring the convergence 
\begin{equation}
 \mathcal{K}_N(A,g) \,\xrightarrow[N\to\infty]{\;\widehat{d}_w\;}\,\overline{\mathcal{K}(A,g)}
\end{equation}
 of the $N$-dimensional truncated Krylov subspace $\mathcal{K}_N(A,g)$ to the (closure of the) full Krylov subspace $\mathcal{K}(A,g)$. Such quite natural \emph{inner Krylov approximability} is actually false when the standard Hausdorff distance is used for an infinite-dimensional $\cH$ \cite[Lemma 7.1]{CM-KrylovPerturbations-2020}. Another natural aspect of the inner Krylov approximability that manifests in the $\widehat{d}_w$-sense is the following \cite[Lemma 7.6]{CM-KrylovPerturbations-2020}:
\begin{equation}
  \left. 
  \begin{array}{c}
   (g_n)_n\,\subset\,\overline{\mathcal{K}(A,g)} \\
   \Vert g_n - g \Vert_\cH \xrightarrow{n \to \infty} 0
  \end{array}
  \right\}\qquad\Rightarrow\qquad   
  \overline{\mathcal{K}(A,g_n)} \,\xrightarrow[n \to \infty]{\;\widehat{d}_w\;}\, \overline{\mathcal{K}(A,g)}\,.
 \end{equation}

 Other natural approximations, instead, do not hold in general in the $\widehat{d}_w$-sense, which suggests that further conditions need be added to the perturbation picture. For example \cite[Lemma 7.3, Examples 7.4 and 7.5]{CM-KrylovPerturbations-2020}, 
 \begin{equation}
    \Vert g_n - g \Vert_\cH \xrightarrow{n \to \infty} 0\qquad\Rightarrow\qquad  
    \begin{cases}
     \;d_w\big( \,\overline{\mathcal{K}(A,g)} \,,\, \overline{\mathcal{K}(A,g_n)}\,\big)\xrightarrow{n \to \infty} 0 \\
     \hspace{1cm}\textrm{but in general}  \\
     \;\overline{\mathcal{K}(A,g_n)} \,\xrightarrow{\quad\widehat{d}_w\qquad} \hspace{-0.8cm}/\hspace{0.7cm} \overline{\mathcal{K}(A,g)}\,.
    \end{cases}
 \end{equation}

 The study of the behaviour of the Krylov solvability along $\widehat{d}_w$-limits is at an early stage. Let us conclude this concise overview with this type of results.

\begin{proposition}\cite[Proposition~7.7]{CM-KrylovPerturbations-2020}\label{prop:Ksolv_perturb}
Let $\cH$ be a separable Hilbert space, and let the following be given:
\begin{itemize}
	\item[(i)] an operator $A \in \mathcal{B}(\cH)$ with inverse $A^{-1} \in \mathcal{B}(\cH)$,
	\item[(ii)] a sequence $(g_m)_{m \in \N} \subset \cH$ such that for each $m \in \N$ the unique solution $f_m := A^{-1}g_m$ to the inverse linear problem $A f_m = g_m$ is a Krylov solution,
	\item[(iii)] a vector $g \in \cH$ such that
	\[\overline{\mathcal{K}(A,g_m)} \xrightarrow{\widehat{d}_w} \overline{\mathcal{K}(A,g)} \,\,\, \text{as } n \to \infty\,.\]
\end{itemize}
Then the unique solution $f:= A^{-1}g$ to the inverse linear problem $Af = g$ is a Krylov solution. If additionally $g_m \to g$ (resp. $g_m \rightharpoonup g$), then $f_m \to f$ (resp. $f_m \rightharpoonup f$).
\end{proposition}

\section{Overview on abstract Friedrichs systems}\label{sec:overviewFriedrichs}


Let us condense in this Section a short overview on Friedrichs systems from their concrete formulation as differential problems to their lifting to abstract inverse linear problems on Hilbert space.

Friedrichs systems were initially introduced by Friedrichs \cite{KOF} in his research on symmetric positive systems, and became relevant in numerical analysis \cite{HMSW,MJensen,Ern-Guermond_book_FiniteElements,EGsemel} as a unified and versatile scheme for solving partial differential equations of different types, including elliptic, parabolic, and hyperbolic.


They are currently being investigated in a variety of directions, including well-posedness \cite{EGC,ABcpde,ABCE,Antonic-Erceg-Mich-2017}, homogenisation \cite{Burazin-Vrdoljak-2014,Burazin-Erceg-Waurick-2023}, representations of boundary conditions \cite{ABcpde}, connection with the classical theory \cite{ABjde,ABCE,ABVnaRWA,ABVisrn}, applications to various initial or boundary value problems of elliptic, hyperbolic, and parabolic type \cite{ABVjmaa,BEmjom,DLS,EGsemel,MDS}, development of different numerical schemes \cite{BDG,TBT,BEF,EGsemel,EGbis,EGter}, classification of abstract realisations of Friedrichs systems \cite{Antonic-Erceg-Mich-2017,EM-FriedrichsDelta2017,Erceg-Soni-2022,Erceg-Soni-2023,Erceg-Soni-2024}, the Friedrichs system structure of contact interactions \cite{EM-FriedrichsDelta2017}, and many others.

The classical Friedrichs operator on a bounded domain $\Omega\subset\mathbb{R}^d$ with Lipschitz boundary $\partial\Omega$ is the first-order differential operator $A$ acting as the $L^2(\Omega,\mathbb{C}^r)\to\mathcal{D}'(\Omega)^r$ map
\begin{equation}\label{eq:classicalFR-1}
 Af\,=\,\sum_{k=1}^d \partial_k (\mathbf{B}_k f)+\mathbf{C} f
\end{equation}
for given matrix functions $\mathbf{B}_k\in W^{1,\infty}(\Omega,M_r)$, $k\in\{1,\dots,d\}$, and $\mathbf{C}\in L^\infty(\Omega,M_r)$ satisfying, for some $\mu>0$, 
\begin{equation}\label{eq:classicalFR-2}
 \mathbf{B}_k\,=\,\mathbf{B}_k^*\qquad\textrm{and}\qquad \mathbf{C}+\mathbf{C}^*+\sum_{k=1}^d \partial_k \mathbf{B}_k\,\geqslant\,\mu\mathbbm{1}\qquad\textrm{a.e. on }\,\Omega\,.
\end{equation}
The inequality in \eqref{eq:classicalFR-2} is in the sense of expectations on $L^2(\Omega,\mathbb{C}^r)$.

The solvability or also the well-posedness of the inverse problem $Af=g$ induced by the operator \eqref{eq:classicalFR-1}-\eqref{eq:classicalFR-2} and by a given $g\in L^2(\Omega,\mathbb{C}^r)$ depend on the kind of realisations of $A$ that are considered in terms of boundary conditions imposed at $\partial\Omega$.

Many other concrete differential problems may be cast in the form $Af=g$ with $A$ of the type \eqref{eq:classicalFR-1}-\eqref{eq:classicalFR-2} \cite{ABcpde,ABjde,ABVnaRWA,ABVjmaa,BEmjom,EGsemel,EGbis}. In addition, it is possible to reformulate all such concrete problems at an abstract Hilbert space setting: this includes the classification of the distinct operator realisations of $A$ (namely, the counterpart of boundary conditions at $\partial\Omega$), as well as the identification of certain intrinsic abstract conditions ensuring the bijectivity of $A$ and hence the solvability of $Af=g$ \cite{EGC,ABcpde,ABCE,Antonic-Erceg-Mich-2017,EM-FriedrichsDelta2017,Erceg-Soni-2022,Erceg-Soni-2023,Erceg-Soni-2024}.

This way one naturally comes to the definition of a \emph{joint pair of abstract Friedrichs operators} on a complex Hilbert space $\cH$ as a pair $(A_0,\widetilde{A}_0)$ of two closed and densely defined operators acting in $\cH$ and 
satisfying the conditions
\begin{equation}\label{eq:defFriedrichs2}
\begin{split}
 & A_0\,\subset\,(\widetilde{A}_0)^*\,=:A_1\,,\qquad \widetilde{A}_0\,\subset\,(A_0)^*\,=:\widetilde{A}_1\,, \\
 & \textrm{$A_1+\widetilde{A}_1$ is a bounded self-adjoint operator in $\cH$} \\
 & \textrm{with strictly positive bottom}\,.
\end{split}
\end{equation}
(The first condition in \eqref{eq:defFriedrichs2} is in the sense of operator inclusion.)

In the concrete case \eqref{eq:classicalFR-1}-\eqref{eq:classicalFR-2} $A_0$ may be taken as a minimally defined version of the operator $A$ therein, and $\widetilde{A}_0$ is its formal adjoint.

In the abstract setting \eqref{eq:defFriedrichs2} it can be shown (\cite[Section 2.1]{EGC}, \cite[Theorem 7]{ABCE}) that the graph norms associated with $A_0$ and $\widetilde{A}_0$ are equivalent, and as a consequence $A_0$ and $\widetilde{A}_0$ have a common (dense) domain
\begin{equation}
 \mathrm{dom}\, A_0\,=\,\mathrm{dom}\,\widetilde{A}_0\,=:\,\mathcal{W}_0
\end{equation}
and so do $A_1$ and $\widetilde{A}_1$, that is,
\begin{equation}
 \mathrm{dom}\, A_1\,=\,\mathrm{dom}\,\widetilde{A}_1\,=:\,\mathcal{W}_1\qquad (\textrm{with } \mathcal{W}_0\,\subset\,\mathcal{W}_1)\,,
\end{equation}
and moreover, 
\begin{equation}
 A_1+\widetilde{A}_1\,\subset\,\overline{A_0+\widetilde{A}_0}\,.
\end{equation}

Furthermore, between the minimal and the maximal operators defined above one is led to consider intermediate pairs $(A, \widetilde{A})$ with
\begin{equation}
 A_0\,\subset\, A \,\subset\, A_1\,,\qquad \widetilde{A}_0\,\subset\, \widetilde{A} \,\subset\, \widetilde{A}_1,
\end{equation}
equivalently, pairs $(V,\widetilde{V})$ of linear subspaces included between $\mathcal{W}_0$ and $\mathcal{W}_1$, and the restrictions $A:=A_1|_V$, $\widetilde{A}:=\widetilde{A}_1|_{\widetilde{V}}$. Each choice of $V$ and $\widetilde{V}$ corresponds, in this abstract setting, to a selection of certain boundary conditions for a concrete differential operator such as \eqref{eq:classicalFR-1}-\eqref{eq:classicalFR-2}.

In particular, a special status is attributed to those $V$'s that are closed in $(\mathcal{W}_1,\|\cdot\|_{A_1})$ and correspondingly those realisations $A_1|_V$'s that are bijective as $(V,\|\cdot\|_{A_1})\to\cH$ maps. Indeed, such operators that are closed and densely defined in $\cH$, are also $(V,\|\cdot\|_{A_1})\to\cH$ bounded, and their inverse is $\cH\to(V,\|\cdot\|_{A_1})$ bounded, implying that the inverse linear problem induced by each such $A_1|_V$ is well-posed. Incidentally, one can see that if $A:=A_1|_V$ is a \emph{closed bijective realisation} of $A_0$, then its adjoint $A^*=(A_1|_V)^*\equiv \widetilde{A}_1|_{\widetilde{V}}$, with $\widetilde{V}:=\mathrm{dom}\,(A_1|_V)^*$, is a closed bijective realisation of $\widetilde{A}_0$.

The investigation and classification of mutually adjoint pairs $(A,A^*)$ of closed bijective realisations of a Friedrichs operator $A_0$ (and of its partner $\widetilde{A}_0$) has been the object of multiple recent studies \cite{Antonic-Erceg-Mich-2017,EM-FriedrichsDelta2017,Erceg-Soni-2022,Erceg-Soni-2023,Erceg-Soni-2024}. Quite noticeably (\cite[Theorem 13]{Antonic-Erceg-Mich-2017}):
\begin{enumerate}
 \item[(i)] A pair of mutually adjoint closed bijective realisations of $(A_0,\widetilde{A}_0)$ always exists.
 \item[(ii)] If both $\ker A_1$ and $\ker \widetilde{A}_1$ are non-trivial, then there exist uncountably many such pairs. If instead at least one among $A_1$ and $\widetilde{A}_1$ is injective, then the only mutually adjoint pair of closed bijective realisations is $(A_1,\widetilde{A}_0)$ when $\ker A_1=\{0\}$, or $(A_0,\widetilde{A}_1)$ when $\ker\widetilde{A}_1=\{0\}$.
\end{enumerate}

 Part (i) above expresses the fact that any differential operator of Friedrichs type admits at least one set of boundary conditions for which the corresponding inverse problem is well-posed. In abstract terms it provides the existence of a `reference' pair $(A_\mathrm{r},A_\mathrm{r}^*)$ with $A_0\subset A_\mathrm{r}\subset A_1$ where $A_\mathrm{r}$ is closed and bijective.

\section[Perspective on Krylov solvability and Friedrichs]{Perspective on Krylov solvability of Friedrichs-type inverse problems}\label{sec:nextresultsandproblems}

As stated already, the study of the Krylov solvability of abstract inverse problems $Af=g$ in a Hilbert space, governed by operators $A$ of Friedrichs type, provides valuable a priori information on the possibility of approximating numerically the solution to a concrete differential inverse problem within the vast class of Friedrichs systems.

At this level of generality the main question is therefore:

\medskip

\textbf{QUESTION 1:} \emph{under what conditions on the (closed) operator $A$ and the vector $g\in\mathrm{ran}\,A$ does the inverse problem $Af=g$ admit solution(s) $f\in\overline{\mathcal{K}(A,g)}$} ?

\medskip

Obviously, Question 1 has a variety of variants, depending on what features of $A$ are considered, such as its general type (normal, self-adjoint, skew-adjoint, or in some special class), or a priori estimates satisfied by $A$ (including boundedness, unboundedness, semi-boundedness, coercivity, etc.).

It is also relevant to investigate in Question 1 the issue of the uniqueness of the Krylov solution, let alone the issue of convergence rate of the Krylov approximants $f_n\to f$ in one of the abstract algorithms considered. 
 
 Work in progress is showing an increasing understanding of the general Question 1 in concrete toy problems as well as special settings that retain a significant amount of instructiveness. 
 
 In this spirit, consider for a moment a generic Hilbert space $\cH$ and an operator $A$ acting in $\cH$ and $g\in\cH$ such that 
 \begin{equation}\label{eq:conditionsqa}
  \begin{split}
   & \textrm{$A$ is densely defined and closed}, \\
   & g\in\mathrm{ran}\,A\cap\mathrm{dom}\,A^*\,, \\
   & A^*g\in\mathcal{D}^{qa}(A^*A)
  \end{split}
 \end{equation}
 (the latter condition meaning that $A^*g$ is a quasi-analytic vector for $A^*A$, in the sense of Definition \ref{def:bdd_a_qa_class}).
 The assumption that $g\in\mathrm{ran}\,A\cap\mathrm{dom}\,A^*$ ensures that $g=Au$ for some $u\in\mathrm{dom}\,A^*A$, and $A^*g=A^*Au\in\mathrm{ran}\,A^*A$. This leads to consider, instead of the inverse problem $Af=g$, the auxiliary inverse problem
 \begin{equation}\label{eq:auxinvprob}
  A^*A f\,=\,A^*g\,.
 \end{equation}
 Here $A^*A$ is self-adjoint (see, e.g., \cite[Theorem 12.11]{Grubb-DistributionsAndOperators-2009}), and $A^*g\in\mathrm{ran}\,A^*A\cap \mathcal{D}^{qa}(A^*A)$. These are precisely the conditions that allow to deduce, via Theorem \ref{th:SelfSkew-adjoint_Ksolvability}, that there exists a unique solution
 \[
  f\,\in\,\overline{\mathcal{K}(A^*A,A^*g)}
 \]
 to the inverse problem \eqref{eq:auxinvprob}.

 On the other hand, one observes that the above solution $f$ also satisfies $Af=g$. Indeed, from $A^*Af=A^*g$ one sees that $Af-g\in\ker A^*$, whereas obviously $Af-g\in\mathrm{ran}\,A\subset(\ker A^*)^\perp$ because $Af$ and $g$ belongs to $\mathrm{ran}\,A$. From the triviality of $\ker A^*\cap (\ker A^*)^\perp$ one thus deduces that $Af=g$.

 It then remains established that under the present conditions \eqref{eq:conditionsqa} there is a unique solution $f$ to $Af=g$ belonging to the (closure of the) Krylov subspace $\overline{\mathcal{K}(A^*A,A^*g)}$.

 This is \emph{not} the actual Krylov subspace one would like to have in the ultimate conclusion of Krylov solvability, namely $\overline{\mathcal{K}(A,g)}$, and the very condition of quasi-analiticity of $A^*g$ with respect to $A^*A$ differs from the standard quasi-analiticity property $g\in\mathcal{D}^{qa}(A)$. The following meaningful question then arises naturally.

 \medskip
 
 \textbf{QUESTION 2:} \emph{Assume that $(A_0,\widetilde{A}_0)$ is a joint pair of abstract Friedrichs operators on $\cH$ satisfying \eqref{eq:conditionsqa} and hence such that there is a unique solution $f$ to $A_0f=g$ belonging to $\overline{\mathcal{K}(A_0^*A_0,A_0^*g)}$. What assumptions need be added in order to obtain the standard Krylov solvability of $Af=g$, where $(A,A^*)$ is one of the closed realisations of $(A_0,\widetilde{A}_0)$ }?
 
 \medskip
 
 In particular, it is of interest to consider Question 2 for the \emph{maximal} realisation $(A_1,\widetilde{A}_1)$ of $(A_0,\widetilde{A}_0)$, in the sense of Section \ref{sec:overviewFriedrichs}, or also for one of the mutually adjoint pairs of closed \emph{bijective} realisations of $(A_0,\widetilde{A}_0)$.

 Current work in progress in collaboration with M.~Erceg (Zagreb) \cite{CEM-2024-KryFri} has allowed us to produce instructive partial answers to Question 2. For instance, quite noticeably, if $A_1+\widetilde{A}_1$ is a multiple of the identity, and if the smooth vectors of $A$ are also smooth vectors of $A^*$ (or vice versa), then indeed
 \[
  \overline{\mathcal{K}(A^*A,A^*g)}\,\subset\,\overline{\mathcal{K}(A,g)}\,,
 \]
 thereby closing in the affirmative the issue of the Krylov solvability for this class of Friedrichs systems.

 This has naturally brought attention to those joint pairs $(A_0,\widetilde{A}_0)$ of abstract Friedrichs operators such that
 \begin{equation}
  A_0 + \widetilde{A}_0 \, = \, \alpha \mathbbm{1} \qquad \alpha\in\mathbb{R}\,.
 \end{equation}
 Consider for instance this very prototypical Friedrichs operator: the complex Hilbert space $\cH=L^2(\mathbb{R})$, and
 \begin{equation}\label{eq:ADp1}
  \mathrm{dom}\,A\,:=\,H^1(\mathbb{R})\,,\qquad Af\,:=\,-f'+f\,.
 \end{equation}
 $A$ is actually an unbounded Friedrichs operator, as it is closed and its adjoint acts as $\frac{\ud}{\ud x}+\mathbbm{1}$ on $H^1(\mathbb{R})$, thereby implying that $A+A^*=2\mathbbm{1}$ is obviously bounded, self-adjoint, and coercive (see    \eqref{eq:defFriedrichs2} above).

 More generally, one may consider variants where $A$ acts as 
  \begin{equation}\label{eq:ADp1-2}
  \begin{split}
     & A\,=\,-\frac{\ud}{\ud x} + c(x) \\
     & \textrm{for $c:\mathbb{R}\to\mathbb{C}$ such that $c(x)+\overline{c(x)}$ is a } \\
     & \textrm{positive bounded function separated from zero}
  \end{split}
 \end{equation}
 (The second term in $-\frac{\ud}{\ud x} + c(x)$ indicates the multiplication by the function $c$.)
 
 \medskip
 
 \textbf{QUESTION 3:} \emph{To investigate the Krylov-solvability of the inverse problem $Af=g$ induced by \eqref{eq:ADp1} or \eqref{eq:ADp1-2} (or variants) under suitable conditions on $g$ and $c$}.
 
 \medskip
 
 It is not difficult to argue that the general arguments preceding Question 2, as well as the above-mentioned partial answers to Question 2, may be made applicable to Question 3. 
 
 For instance, we could show \cite{CEM-2024-KryFri} that in the case \eqref{eq:ADp1} it suffices to assume further that the $\widehat{g}\in C^\infty_c(\mathbb{R})$ (i.e., the Fourier transform of $g$ is smooth and compactly supported) in order to deduce that $A^*g$ is a \emph{bounded} (and hence quasi-analytic) vector for $A^*A$, in the sense of Definition \ref{def:bdd_a_qa_class}. This ensures that $A^*g$ is in particular quasi-analytic for $A^*A$ and as a consequence of the facts established before Question 2 and the partial results stated thereafter one deduces that the Friedrichs-type inverse problem
 \[
  -f'+f\,=\,g
 \]
 admits a unique Krylov solution. Analogous partial result are currently being produced for the variant \eqref{eq:ADp1-2}.

 In fact, the latter line of reasoning mirrors at the level of concrete differential operators what Questions 1 and 2 deal with at an abstract setting, and could be therefore continued by inquiring the following.
 
  \medskip
 
 \textbf{QUESTION 4:} \emph{To investigate the Krylov-solvability of the inverse problem $Af=g$ induced by the classical Friedrichs operator \eqref{eq:classicalFR-1} under suitable conditions on the matrix functions  matrix functions $\mathbf{B}_k\in W^{1,\infty}(\Omega,M_r)$, $k\in\{1,\dots,d\}$, and $\mathbf{C}\in L^\infty(\Omega,M_r)$}.
 
 \medskip
 
 Back to the abstract Hilbert space setting, a variety of additional meaningful problems arise naturally. Among them, we find the following two particularly intriguing and promising.

   \medskip
 
 \textbf{QUESTION 5:} \emph{To characterise the Krylov intersection $\mathcal{I}(A,g)$ relative to an abstract Friedrichs operator $A$ acting in a Hilbert space $\cH$ and a vector $g\in\cH$ (see \eqref{eq:def_Kint} and \eqref{eq:def_KInt_unbdd} above), as well as the Krylov reducibility of $A$ (as in the spirit of Proposition \ref{prop:Kreduced_Ksolvable}). Or, at the concrete level, to characterise the Krylov intersection and the Krylov reducibility of the classical Friedrichs operator \eqref{eq:classicalFR-1}}.
 
    \medskip
    
 Question 5 is a `structural' question, where the general goal is to investigate when the Krylov intersection is trivial, which under suitable conditions is a sufficient criterion for Krylov solvability, in the spirit of Propositions \ref{prop:trivialKint_Ksolv} above. It is reasonable to expect that the special constraints \eqref{eq:defFriedrichs2} in the abstract Hilbert space setting (or \eqref{eq:classicalFR-1}-\eqref{eq:classicalFR-2} in a concrete differential setting) on an operator of Friedrichs type result in a special behaviour concerning its Krylov reducibility or the triviality of its Krylov intersection, which in turn determine the Krylov solvability of the associated inverse problem.

    \medskip
 
 \textbf{QUESTION 6:} \emph{Consider a perturbation $\widetilde{A}\widetilde{f}=\widetilde{g}$ of the inverse linear problem $Af=g$ induced by a Friedrichs operator $A$, where in some suitable sense $\widetilde{A}\sim A$ and $\widetilde{g}\sim g$. When, in the $\widehat{d}_w$-sense (Section \ref{subsec:KrylovPerturb}), are the Krylov subspaces $\mathcal{K}(A,g)$ and $\mathcal{K}(\widetilde{A},\widetilde{g})$ close ? When does the perturbation preserve or alter the Krylov solvability of the problem ? What is the answer to the above questions when, in particular, the perturbation preserves the Friedrichs character of the linear operator ?}
 
    \medskip
 
 Concerning Question 6, it was recalled in Section \ref{subsec:KrylovPerturb} that the sole vicinity $\widetilde{A}\sim A$ or $\widetilde{g}\sim g$, even in some strong sense, does not guarantee the smallness of 
 \[
  \widehat{d}_w \Big( \,\overline{\mathcal{K}(A,g)} \,,\, \overline{\mathcal{K}(\widetilde{A},\widetilde{g})}\,\Big)\,,
 \]
 and, typically, it is the term $d_w \big( \, \overline{\mathcal{K}(\widetilde{A},\widetilde{g})}\,,\, \overline{\mathcal{K}(A,g)}\,\big)$ in \eqref{eq:defweakgap} to be problematic. It was also recalled (Proposition \ref{prop:Ksolv_perturb}) that if the perturbation leaves the Krylov subspaces $\overline{\mathcal{K}(A,g)}$ and $\overline{\mathcal{K}(\widetilde{A},\widetilde{g})}$ $\widehat{d}_w$-close, then under favourable conditions the Krylov solvability of the corresponding inverse problems is preserved, whereas in general it may be gained or lost in a perturbation (Example \ref{eg:Krylovperturb_lossgain}). It is plausible to expect, once again, that the constraints for $A$ to be an operator of Friedrichs type select distinguished results in this general scenario. It should be even more so for those perturbations that actually \emph{preserve} the status of Friedrichs system of the inverse problem.

 The above selection of questions, be they formulated at an abstract or concrete level, should convey the message, as we strongly believe, that this intriguing subject is ready for very promising developments. For one more time we should like to stress that general criteria of Krylov solvability in the context of inverse problems of Friedrichs type are most valuable, in view of the wide popularity of Krylov approximation methods in numerical analysis, as well as of the treatment of partial differential equations cast in the form of Friedrichs systems.
 
\section{Acknowledgements}
This work is partially supported by the Italian National Institute for Higher Mathematics INdAM and the Alexander von Humboldt Foundation, Bonn.


\begin{thebibliography}{10}

\bibitem{ABcpde}
{\sc N.~Antoni{\'c} and K.~s. Burazin}, {\em {Intrinsic boundary conditions for
  {F}riedrichs systems}}, Comm. Partial Differential Equations, 35 (2010),
  pp.~1690--1715.

\bibitem{ABjde}
\leavevmode\vrule height 2pt depth -1.6pt width 23pt, {\em {Boundary operator
  from matrix field formulation of boundary conditions for {F}riedrichs
  systems}}, J. Differential Equations, 250 (2011), pp.~3630--3651.

\bibitem{ABCE}
{\sc N.~Antoni{\'c}, K.~s. Burazin, I.~Crnjac, and M.~Erceg}, {\em {Complex
  {F}riedrichs systems and applications}}, J. Math. Phys., 58 (2017),
  pp.~101508, 22.

\bibitem{ABVisrn}
{\sc N.~Antoni{\'c}, K.~s. Burazin, and M.~Vrdoljak}, {\em {Connecting
  classical and abstract theory of {F}riedrichs systems via trace operator}},
  ISRN Math. Anal.,  (2011), pp.~Art. ID 469795, 14.

\bibitem{ABVjmaa}
\leavevmode\vrule height 2pt depth -1.6pt width 23pt, {\em {Heat equation as a
  {F}riedrichs system}}, J. Math. Anal. Appl., 404 (2013), pp.~537--553.

\bibitem{ABVnaRWA}
\leavevmode\vrule height 2pt depth -1.6pt width 23pt, {\em {Second-order
  equations as {F}riedrichs systems}}, Nonlinear Anal. Real World Appl., 15
  (2014), pp.~290--305.

\bibitem{Antonic-Erceg-Mich-2017}
{\sc N.~Antoni{\'c}, M.~Erceg, and A.~Michelangeli}, {\em {Friedrichs systems
  in a {H}ilbert space framework: {S}olvability and multiplicity}}, J.
  Differential Equations, 263 (2017), pp.~8264--8294.


\bibitem{Atkinson-Han-TheoNumAnal2009}
{\sc K.~Atkinson and W.~Han}, {\em {Theoretical numerical analysis}}, vol.~39
  of {Texts in Applied Mathematics}, Springer, Dordrecht, third~ed., 2009.
\newblock A functional analysis framework.

\bibitem{TBT}
{\sc T.~Bui-Thanh}, {\em {From {G}odunov to a unified hybridized discontinuous
  {G}alerkin framework for partial differential equations}}, J. Comput. Phys.,
  295 (2015), pp.~114--146.

\bibitem{BDG}
{\sc T.~Bui-Thanh, L.~Demkowicz, and O.~Ghattas}, {\em {A unified discontinuous
  {P}etrov-{G}alerkin method and its analysis for {F}riedrichs' systems}}, SIAM
  J. Numer. Anal., 51 (2013), pp.~1933--1958.

\bibitem{Burazin-Erceg-Waurick-2023}
{\sc K.~Burazin, M.~Erceg, and M.~Waurick}, {\em {G-convergence of Friedrichs
  systems revisited}}, arXiv:2307.01552 (2023).

\bibitem{Burazin-Vrdoljak-2014}
{\sc K.~Burazin and M.~Vrdoljak}, {\em {Homogenisation theory for {F}riedrichs
  systems}}, Commun. Pure Appl. Anal., 13 (2014), pp.~1017--1044.

\bibitem{BEmjom}
{\sc K.~s. Burazin and M.~Erceg}, {\em {Non-stationary abstract {F}riedrichs
  systems}}, Mediterr. J. Math., 13 (2016), pp.~3777--3796.

\bibitem{BEF}
{\sc E.~Burman, A.~Ern, and M.~A. Fern{\'a}ndez}, {\em {Explicit
  {R}unge-{K}utta schemes and finite elements with symmetric stabilization for
  first-order linear {PDE} systems}}, SIAM J. Numer. Anal., 48 (2010),
  pp.~2019--2042.

\bibitem{CEM-2024-KryFri}
{\sc N.~A. Caruso, M.~Erceg, and A.~Michelangeli}, {\em {Work in progress}},
  (2024).

\bibitem{CM-2019_ubddKrylov}
{\sc N.~A. Caruso and A.~Michelangeli}, {\em {Krylov {S}olvability of
  {U}nbounded {I}nverse {L}inear {P}roblems}}, Integral Equations Operator
  Theory, 93 (2021), p.~Paper No. 1.

\bibitem{CM-KrylovMonograph-2021}
\leavevmode\vrule height 2pt depth -1.6pt width 23pt, {\em {Inverse linear
  problems on {H}ilbert space and their {K}rylov solvability}}, {Springer
  Monographs in Mathematics}, Springer, Cham, [2021] \copyright 2021.

\bibitem{CM-Nemi-unbdd-2019}
\leavevmode\vrule height 2pt depth -1.6pt width 23pt, {\em {Convergence of the
  conjugate gradient method with unbounded operators}}, Operators and Matrices,
  16 (2022), pp.~35--68.

\bibitem{CM-KrylovPerturbations-2020}
\leavevmode\vrule height 2pt depth -1.6pt width 23pt, {\em {Krylov solvability
  under perturbations of abstract inverse linear problems}}, Journal of Applied Analysis, 29 (2023), pp.~3--29.

\bibitem{CMN-2018_Krylov-solvability-bdd}
{\sc N.~A. Caruso, A.~Michelangeli, and P.~Novati}, {\em {On Krylov solutions
  to infinite-dimensional inverse linear problems}}, Calcolo, 56 (2019), p.~32.

\bibitem{CMN-truncation-2018}
\leavevmode\vrule height 2pt depth -1.6pt width 23pt, {\em {On general
  projection methods and convergence behaviours for abstract linear inverse
  problems}}, Asymptotic Analysis, 127 (2022), pp.~167--189.

\bibitem{C-KrylovNormal-2022}
{\sc {Caruso, No{\`e} Angelo}}, {\em {A note on the Krylov solvability of
  compact normal operators on Hilbert space}}, {Complex Analysis and Operator
  Theory}, 17 (2023).

\bibitem{Chen3-1997}
{\sc M.~Chen, Z.~Chen, and G.~Chen}, {\em {Approximate solutions of operator
  equations}}, vol.~9 of {Series in Approximations and Decompositions}, World
  Scientific Publishing Co., Inc., River Edge, NJ, 1997.

\bibitem{Cipra-SIAM-News}
{\sc B.~A. Cipra}, {\em {The best of the 20th century: Editors name top 10
  algorithms}}, SIAM News, 33 (2005).

\bibitem{DLS}
{\sc B.~Despr{\'e}s, F.~Lagouti{\`e}re, and N.~Seguin}, {\em {Weak solutions to
  {F}riedrichs systems with convex constraints}}, Nonlinearity, 24 (2011),
  pp.~3055--3081.

\bibitem{Dongarra-Sullivan-Best10-2000}
{\sc J.~Dongarra and F.~Sullivan}, {\em {The Top 10 Algorithms (Guest editors'
  intruduction)}}, Comput. Sci. Eng., 2 (2000), pp.~22--23.

\bibitem{Engl-Hanke-Neubauer-1996}
{\sc H.~W. Engl, M.~Hanke, and A.~Neubauer}, {\em {Regularization of inverse
  problems}}, vol.~375 of {Mathematics and its Applications}, Kluwer Academic
  Publishers Group, Dordrecht, 1996.

\bibitem{EM-FriedrichsDelta2017}
{\sc M.~Erceg and A.~Michelangeli}, {\em {On contact interactions realised as
  {F}riedrichs systems}}, Complex Anal. Oper. Theory, 13 (2019), pp.~703--736.

\bibitem{Erceg-Soni-2022}
{\sc M.~Erceg and S.~K. Soni}, {\em {Classification of classical {F}riedrichs
  differential operators: one-dimensional scalar case}}, Commun. Pure Appl.
  Anal., 21 (2022), pp.~3499--3527.

\bibitem{Erceg-Soni-2023}
\leavevmode\vrule height 2pt depth -1.6pt width 23pt, {\em {The von Neumann
  extension theory for abstract Friedrichs operators }}, arXiv:2312.09618
  (2023).

\bibitem{Erceg-Soni-2024}
\leavevmode\vrule height 2pt depth -1.6pt width 23pt, {\em {Friedrichs systems
  on an interval }}, arXiv:2401.11941 (2024).

\bibitem{Ern-Guermond_book_FiniteElements}
{\sc A.~Ern and J.-L. Guermond}, {\em {Theory and practice of finite
  elements}}, vol.~159 of {Applied Mathematical Sciences}, Springer-Verlag, New
  York, 2004.

\bibitem{EGsemel}
{\sc A.~Ern and J.-L. Guermond}, {\em {Discontinuous {G}alerkin methods for
  {F}riedrichs' systems. {I}. {G}eneral theory}}, SIAM J. Numer. Anal., 44
  (2006), pp.~753--778.

\bibitem{EGbis}
{\sc A.~Ern and J.-L. Guermond}, {\em {Discontinuous {G}alerkin methods for
  {F}riedrichs' systems. {II}. {S}econd-order elliptic {PDE}s}}, SIAM J. Numer.
  Anal., 44 (2006), pp.~2363--2388.

\bibitem{EGter}
\leavevmode\vrule height 2pt depth -1.6pt width 23pt, {\em {Discontinuous
  {G}alerkin methods for {F}riedrichs' systems. {III}. {M}ultifield theories
  with partial coercivity}}, SIAM J. Numer. Anal., 46 (2008), pp.~776--804.

\bibitem{EGC}
{\sc A.~Ern, J.-L. Guermond, and G.~Caplain}, {\em {An intrinsic criterion for
  the bijectivity of {H}ilbert operators related to {F}riedrichs' systems}},
  Comm. Partial Differential Equations, 32 (2007), pp.~317--341.

\bibitem{KOF}
{\sc K.~O. Friedrichs}, {\em {Symmetric positive linear differential
  equations}}, Comm. Pure Appl. Math., 11 (1958), pp.~333--418.

\bibitem{Grubb-DistributionsAndOperators-2009}
{\sc G.~Grubb}, {\em {Distributions and operators}}, vol.~252 of {Graduate
  Texts in Mathematics}, Springer, New York, 2009.

\bibitem{Hanke-ConjGrad-1995}
{\sc M.~Hanke}, {\em {Conjugate gradient type methods for ill-posed problems}},
  vol.~327 of {Pitman Research Notes in Mathematics Series}, Longman Scientific
  \& Technical, Harlow, 1995.

\bibitem{Hansen-Illposed-1998}
{\sc P.~C. Hansen}, {\em {Rank-deficient and discrete ill-posed problems}},
  {SIAM Monographs on Mathematical Modeling and Computation}, Society for
  Industrial and Applied Mathematics (SIAM), Philadelphia, PA, 1998.
\newblock Numerical aspects of linear inversion.

\bibitem{HMSW}
{\sc P.~Houston, J.~A. Mackenzie, E.~S{\"u}li, and G.~Warnecke}, {\em {A
  posteriori error analysis for numerical approximations of {F}riedrichs
  systems}}, Numer. Math., 82 (1999), pp.~433--470.

\bibitem{MJensen}
{\sc M.~Jensen}, {\em {Discontinuous Galerkin methods for Friedrichs systems
  with irregular solutions}}, Ph.D.~thesis, University of Oxford (2004).

\bibitem{Kato-perturbation}
{\sc T.~Kato}, {\em {Perturbation theory for linear operators}}, {Classics in
  Mathematics}, Springer-Verlag, Berlin, 1995.
\newblock Reprint of the 1980 edition.

\bibitem{kranoselski-1972-approxsoll}
{\sc M.~A. {Krasnosel\cprime ski\u{\i}}, G.~M. Va\u{\i}nikko, P.~P.
  Zabre\u{\i}ko, Y.~B. Rutitskii, and V.~Y. Stetsenko}, {\em {Approximate
  solution of operator equations}}, Wolters-Noordhoff Publishing, Groningen,
  1972.
\newblock Translated from the Russian by D. Louvish.

\bibitem{Liesen-Strakos-2003}
{\sc J.~Liesen and Z.~e. {Strako\v s}}, {\em {Krylov subspace methods}},
  {Numerical Mathematics and Scientific Computation}, Oxford University Press,
  Oxford, 2013.
\newblock Principles and analysis.

\bibitem{MDS}
{\sc C.~Mifsud, B.~Despr{\'e}s, and N.~Seguin}, {\em {Dissipative formulation
  of initial boundary value problems for {F}riedrichs' systems}}, Comm. Partial
  Differential Equations, 41 (2016), pp.~51--78.

\bibitem{Nemirovskiy-Polyak-1985}
{\sc A.~S. Nemirovskiy and B.~T. Polyak}, {\em {Iterative methods for solving
  linear ill-posed problems under precise information. {I}}}, Izv. Akad. Nauk
  SSSR Tekhn. Kibernet.,  (1984), pp.~13--25, 203.

\bibitem{Nemirovskiy-Polyak-1985-II}
\leavevmode\vrule height 2pt depth -1.6pt width 23pt, {\em {Iterative methods
  for solving linear ill-posed problems under precise information. {II}}},
  Engineering Cybernetics, 22 (1984), pp.~50--57.

\bibitem{Nevanlinna-1993-converg-iterat-book}
{\sc O.~Nevanlinna}, {\em {Convergence of iterations for linear equations}},
  {Lectures in Mathematics ETH Z{\"u}rich}, Birkh{\"a}user Verlag, Basel, 1993.

\bibitem{Quarteroni-book_NumModelsDiffProb}
{\sc A.~Quarteroni}, {\em {Numerical models for differential problems}},
  vol.~16 of {MS\&A. Modeling, Simulation and Applications}, Springer, Cham,
  2017.
\newblock Third edition.

\bibitem{Saad-2003_IterativeMethods}
{\sc Y.~Saad}, {\em {Iterative methods for sparse linear systems}}, Society for
  Industrial and Applied Mathematics, Philadelphia, PA, second~ed., 2003.

\bibitem{schmu_unbdd_sa}
{\sc K.~Schm{\"u}dgen}, {\em {Unbounded self-adjoint operators on {H}ilbert
  space}}, vol.~265 of {Graduate Texts in Mathematics}, Springer, Dordrecht,
  2012.

\bibitem{Sullivan-UncertQuant}
{\sc T.~J. Sullivan}, {\em {Introduction to uncertainty quantification}},
  vol.~63 of {Texts in Applied Mathematics}, Springer, Cham, 2015.

\bibitem{vainikko1981}
{\sc G.~Vainikko}, {\em {Regular convergence of operators and approximate
  solution of equations}}, Journal of Soviet Mathematics, 15 (1981),
  pp.~675--705.

\bibitem{Vainikko-1993}
\leavevmode\vrule height 2pt depth -1.6pt width 23pt, {\em {Multidimensional
  weakly singular integral equations}}, vol.~1549 of {Lecture Notes in
  Mathematics}, Springer-Verlag, Berlin, 1993.

\bibitem{vainikko1974}
{\sc G.~Vainikko and O.~Karma}, {\em {The convergence of approximate methods
  for solving linear and non-linear operator equations}}, USSR Computational
  Mathematics and Mathematical Physics, 14 (1974), pp.~9--19.

\end{thebibliography}
\def\cprime{$'$}

\end{document}